\documentclass[12pt]{amsart}
\usepackage{amsmath,amssymb}
\usepackage{amsmath,amssymb,enumerate}
\usepackage{epsf,epsfig,amsfonts,graphicx,color,a4wide}  
 \usepackage[applemac]{inputenc} 

 \date{}   
  
\newtheorem{theorem}{\rm\bf Theorem}
\newtheorem*{spec}{\rm\bf Special Case}

\newtheorem*{fact}{\rm\bf Fact}

\theoremstyle{definition}

\newtheorem*{Rem}{\rm\bf Remark}

\newcommand{\weg}[1]{}

\begin{document} 

\title[Locally symmetric Finsler spaces]{There exist no locally symmetric Finsler spaces of positive or negative flag  curvature}

\author{Vladimir S. Matveev} 
\address{Mathematisches Institut, Friedrich-Schiller Universit\"at Jena\\
07737 Jena, Germany}  
\email{vladimir.matveev@uni-jena.de}

\begin{abstract}  We show that the results of  Foulon \cite{foulon,foulon2} and Kim \cite{kim}
(independently, Deng and Hou \cite{deng})  about the nonexistence of locally symmetric Finsler metrics of  positive or negative flag  curvature are in fact local. \end{abstract} 
\maketitle

  \emph{Finsler metric}  on a smooth manifold  $M$  is   a  continuous 
function $F:TM\to [0,\infty)$ such that for every point $p\in M$ the  restriction $F_p = F_{|T_pM}$ on the tangent space at $p$ is a  
{Minkowski norm}, that is $F_p$  is positively homogenous and convex and it vanishes only at $y = 0$: 
\begin{enumerate}[ \ (a)]
  \item $F_p(\lambda \cdot y) = \lambda \cdot   F_p (y) $ for any $\lambda \geq 0$.
   \item $F_p (y+ \tilde y ) \le F_{p} (y) + F_{p} (\tilde y)$.
   \item $F_p (y)= 0 $ \  $ \Rightarrow$ \  $y=0$.  
  \end{enumerate}  

We we will also  assume that our   Finsler metric   is  \emph{of class $C^2$}, that is   the restriction of $F$ to the slit tangent bundle $TM^0 = TM\setminus (\text{the zero section})$ is a function of class  $C^2$, and that  it is   \emph{strongly convex},  that is the Hessian of the restriction of $\tfrac{1}{2}F^2$ to $T_pM\setminus \{0 \}$  (which will be denoted by $g= g_{ij}$ later) 
is   positively  definite for any 
$p \in M$  and any nonzero vector $y\in T_pM$.

The Finsler manifold $(M,F)$ is called \emph{locally symmetric},   if for every point $p\in M$ there exists $r=r(p)>0$ and an isometry $ I_p: B_r(p)\to B_r(p)$ (called the \emph{reflection} at $p$) such that $ I_p(p) = p$ and $d_p( I_p) = -\mathrm{id} : T_pM \to T_pM$. Here $B_r(p)$ denotes the ball of radius $r$ around $p$.

Our main result is:
\begin{theorem} \label{thm}
Suppose the flag curvature of a locally symmetric Finsler 
 metric is negative or positive. Then,  the metric is actually a Riemannian metric, that is,  there exists a Riemannian metric $h$ such that for all $p\in M$, $y \in T_pM$ we have $F_p(y)= \sqrt{h_p(y, y) }$. 
\end{theorem} 

Note that there exist examples of locally (and even globally) symmetric Finsler metrics such that the flag curvature 
changes the sign, or is   nonpositive, or is nonnegative. 
Actually, the (reversible) Minkowski space is  already  an example of globally symmetric Finsler space of nonpositive and of nonnegative flag curvature, since its flag curvature is  zero. One can also take the direct product of a Minkowski space with the round sphere and/or with the hyperbolic space 
 constructing  examples such that the flag curvature is nonpositive and is negative somewhere, or such that the flag curvature  is nonnegative and is positive  somewhere, or such that it changes the sign. The holonomy group of all these examples is reducible but one can also construct irreducible examples by perturbing a locally symmetric  Riemannian   metric of rank $\ge 2 $  by an arbitrary (smooth,  small and homogeneous of degree 1) function of the Chevalley's polynomials. Since the  Chevalley's polynomials are preserved by  any isometry of the initial Riemannian metric
(see e.g. \cite{palais,Ch}),  the obtained  Finsler metric is still locally symmetric. 

Special cases of  Theorem \ref{thm}  are  the following two statements  which were known before in the ``global'' setting;   the initial proof of these statements is also  ``global'', i.e., it   requires the  assumption that the  manifold $M$ is compact,  and is very different from our proof which is local. 

\begin{spec}[Corollary 1  of Foulon \cite{foulon2};  also follows  from theorem (A) of Foulon \cite{foulon}]  Let $(M,F)$ be a Finsler locally symmetric space of negative flag curvature. If $M$ is compact, then  $F$ is actually a Riemannian metric.\end{spec}

\begin{spec}[independently Deng and Hou \cite{deng} and Kim \cite{kim}]  Let $(M,F)$ be a Finsler locally symmetric space of positive flag curvature. If $M$ is compact, then  $F$ is actually a Riemannian metric.\end{spec}

\begin{Rem}  Actually, Foulon in \cite{foulon} has a less restrictive definition of locally symmetric spaces (locally symmetric spaces in the our definition are also locally symmetric in the definition of \cite{foulon} but not vice versa) so his result is in fact stronger and we can not repeat it by our methods or prove its local version.
\end{Rem} 

\subsubsection*{ Proof of Theorem \ref{thm}. } Our proof of Theorem \ref{thm} is  based on the following  recent result: 

\begin{fact}[\cite{troyanov}, Remark (A) in \S 8 + Theorem 9.2] Let $(M, F)$ be a $C^2$-smooth
Finsler manifold. If $(M, F)$ is locally
symmetric, then $F$ is Berwald.  Moreover, the associated connection is the Levi-Civita connection of a locally symmetric Riemannian metric. (We denote this Riemannian metric by   $h$, in paper \cite{troyanov} it is called the Binet-Legendre metric associated to the Finsler metric.) \end{fact} 

Recall that a Finsler manifold is \emph{Berwald}, if there exists   a torsion free linear connection $\nabla$  called \emph{associated connection}   such that the parallel transport preserves  the Finsler metric.

It is well-known  (see for example \cite{Bao,bao1})
that the flag curvature $K_p(y,V)$ for the Berwald metrics can be calculated by the following procedure. Compute the curvature tensor $R^{j}_{\ ik\ell}$
 of the associated  connection;  for every $p$ we view the curvature tensor 
 as the mapping \begin{equation} \label{2} R: T_pM\times T_pM \times T_pM  \to 
 T_pM, \ R(a,b)c= c^iR^j_{\ i k \ell} a^k b^\ell. \end{equation} 
  Then, for every two linearly independent 
  vectors $y^i, V^i\in T_pM$  we have 
  \begin{equation} \label{flag} K_p(y,V)=    \frac{  g_y (V,  R (V , y)y) }{ g_y(y,y)g_y(V,V)- g_y(V,y)^2}.\end{equation}     
 Here $g_y$  is the second differential of the restriction of 
the  function $\tfrac{1}{2}F^2$ to $T_pM$; 
 $g_y= g_{ij} $ is a $(0,2)$-tensor whose components  depend on the point $p$ in $M$  and on the tangent vector  $y\in T_pM$.

 Let us now consider the Riemannian metric  $h$ whose existence we recalled in Fact above: it is  locally symmetric and 
 its  Levi-Civita connection is the associated connection of $F$. Let us now 
 show that if the flag curvature of $F$ is positive (for all linearly independent $y$ and $V\in T_pM $)  
 then the sectional curvature of $h$ is also  positive, 
 and   if the flag curvature of $F$ is negative 
 then the sectional curvature of $h$ is also negative.

In order to do this,  for each vector $y\in T_xM$, let us consider the  endomorphism 
 $$
A_y:T_pM\to T_pM, \  V\mapsto  R (V ,y )y.
 $$
 
In the tensor notation, $(A_y)^j_k= R^j_{\ i k \ell } y^iy^\ell$. 
  Since the bilinear form $(\xi,\nu)\mapsto    h(\xi, R ( \nu, y)y)$ is symmetric with  respect to $\xi$ and $\nu$ 
   because of the symmetries of the curvature tensor, 
   for each $y$  the endomorphism $A_y$   is diagonalizable. Clearly, $y$ is an eigenvector 
    of $A_y$ with eigenvalue $0$. 
   Comparing the formula for $A_y$  with \eqref{flag}, 
 we see that the flag curvature is given by 
 $$
 K_p(y,V)= \frac{g_y(V, A_y(V)) }{ g_y(y,y)g_y(V,V)- g_y(V,y)^2}. 
$$ 

If $K_p(y,V)$ is  positive for all linearly independent $y$ and $V$,
 then for each $y\ne 0$ 
all eigenvalues of $A_y$ except of the eigenvalue $0$ corresponding to the eigenvector $y$  are positive, otherwise  the pair   $(y, V)$ where $V$ is an eigenvector 
with nonpositive eigenvalue has $K_p(y, V)\le 0$. 
  Then, the sectional curvature of $h$ (which is given by 
$\frac{h(V, A_y(V)) }{ h(y,y)h(V,V)- h(V,y)^2}$) 
is also positive  for all linearly independent $y$ and $V$ as we  claim.

The case when $K_p(y,V)$ is negative for all linearly independent $y$ and $V$  is virtually the same -- one needs to replace ``positive'' by ``negative'' in the previous arguments.

Finally, if the flag curvature is positive, or if it is negative, the sectional curvature of $h$  is positive resp. negative, i.e., is never equal to zero. 
But then the holonomy group of $h$ is transitive by \cite[Theorem 9]{simons}. 
Since both 
the Finsler function 
$F$ and the metric $h$  are  preserved by by the holonomy group, the function $\frac{F_p(y)^2}{h_p(y,y)}  $ does not depend on $y$ so the metric $F$ is a Riemannian metric as we claimed. 
\emph{    Theorem \ref{thm} is proved}.

\subsubsection*{  Acknowledgments.} The paper was started during the conference 
 XVII Escola de Geometria Diferencial at Manaus and was initiated by  discussions with J.-H. Eschenburg. The last part of the proof of Theorem \ref{thm}  was explained to me  by S. Deng during my visit to Nankai University. I am grateful to the organisational commettee of the above mentioned conference and  to Nankai University  for partial financial support and their hospitality. I also thank  Friedrich-Schiller    University of Jena  for partial financial support.


\end{document}